\documentclass[12pt]{amsart}
\usepackage{epsfig}
\newfont{\sheaf}{eusm10 scaled\magstep1}

\newcommand{\ra}{\ensuremath{\rightarrow}}

\def\eea{\end{eqnarray*}}
\def\bea{\begin{eqnarray*}}

\def\X{{\mathcal{X}}}
\def\ZZZ{{\mathcal{Z}}}
\def\E{{\mathcal{E}}}
\def\Y{{\mathcal{Y}}}
\def\K{{\mathcal{K}}}
\def\MM{{\mathcal{M}}}
\def\de{{\delta}}
\def\De{{\Delta}}

\newcommand{\Proof}{{\it Proof. }}
\newcommand{\QED}{{\hfill $Q.E.D.$}}

\newtheorem{teo}{Theorem}[section]
\newtheorem{df}[teo]{Definition}
\newtheorem{lem}[teo]{Lemma}

\newtheorem{oss}[teo]{Remark}

\newcommand{\C}{\ensuremath{\mathbb{C}}}
\newcommand{\R}{\ensuremath{\mathbb{R}}}
\newcommand{\Z}{\ensuremath{\mathbb{Z}}}
\newcommand{\Q}{\ensuremath{\mathbb{Q}}}

\newcommand{\hol}{\ensuremath{\mathcal{O}}}

\newcommand{\PP}{\ensuremath{\mathbb{P}}}

\newcommand{\D}{\ensuremath{\mathcal{D}}}
\newcommand{\CCC}{\ensuremath{\mathcal{C}}}

\newcommand{\sM}{{\mathcal M}}
\newcommand{\sK}{{\mathcal K}}

\begin{document}

\title[Canonical symplectic structure and deformations]{Canonical symplectic
structures and deformations of  algebraic surfaces.}

\author{Fabrizio Catanese\\
   Universit\"at Bayreuth}

\indent
This article is dedicated to the memory of Boris  Moisezon

\date{May 9, 2007}

\begin{abstract}
    We show that a minimal surface of general type has a canonical
symplectic structure (unique up to symplectomorphism)
   which is  invariant for smooth deformation.
   We show that the symplectomorphism type is also invariant
   for deformations which  allow certain normal singularities, provided one
remains in the same smoothing component.

We use this technique to show that the Manetti surfaces
yield examples of surfaces of general type which are
not deformation equivalent but are
canonically symplectomorphic.

    \end{abstract}
\maketitle

\section{Introduction}

It is well known that if two smooth compact complex manifolds $M, M'$
are deformation equivalent, then there is a diffeomorphism
$ f : M' \ra M$ such that $ f^* (c_1(K_M)) = c_1(K_{M'}) $
(we refer  the reader to \cite{cime} for  a general discussion
of these type of problems and for more complete references).

It was for long time an open question (called Def= Diff question) whether
conversely two diffeomorphic
minimal algebraic surfaces $S, S'$ would be deformation equivalent
(see \cite{katata} and \cite{f-m1}).
   Note that  Seiberg-Witten
theory  furnished a simple proof of the fact, predicted by Donaldson theory
much earlier,
that a diffeomorphism $ f : S' \ra S$ would have the property that
$ f^* (c_1(K_S)) = \pm c_1(K_{S'})$. There are indeed 'easy'
counterexamples to the Def= Diff question, given by pairs of such 
surfaces $S, S'$,
admitting a diffeomorphism $ f : S' \ra S$ with
$ f^* (c_1(K_S)) = - c_1(K_{S'})$, but none with
$ f^* (c_1(K_S)) =  c_1(K_{S'})$ (\cite{k-k}, \cite{cat4}, \cite{bcg}).

The first counterexamples to the refined Def= Diff conjecture where
a diffeomorphism $ f : S' \ra S$ with
$ f^* (c_1(K_S)) =  c_1(K_{S'})$ is required were given by Manetti
in \cite{man4}, and we shall refer to the surfaces constructed
in that paper as {\bf Manetti surfaces}.

Manetti surfaces are  not simply connected,
and one could suspect that their  universal  abelian covers (which are compact)
could be deformation equivalent.

For this reason \cite{cat02} proposed some examples of simply
connected minimal surfaces of general type,
  called 'abc' examples, showing that they are not
deformation equivalent.

Later on, our joint paper with Wajnryb
\cite{c-w} showed that 'abc' surfaces  with fixed integer
  constants $b $ and $a+c$
yield diffeomorphic surfaces of general type.
Thus \cite{c-w} exhibited simply connected counterexamples to
the refined Def= Diff question.

On the other hand, in  \cite{cat02} it was shown that
a minimal surface of general type $S$ has a canonical
symplectic structure (unique up to symplectomorphism)
   which is  invariant for smooth deformation, and which
is called 'canonical' because its
cohomology class  is exactly the class of the canonical divisor $K_S$.

A further refinement of the Def = Diff  question is then
whether   minimal algebraic surfaces of general
type which are canonically symplectomorphic are also deformation equivalent.

We prove in this paper that the Manetti surfaces yield counterexamples
also to this further question.

The same question for the simply connected case  remains open,
and we have not yet been able to decide whether
the 'abc' examples are canonically symplectomorphic.

For this reason we take up here again the Manetti surfaces,
already considered in \cite{cat02}.
Our treatment here is direct and elementary, with the purpose
of being understood both by algebraic geometers and
by symplectic geometers.
First of all we  avoid to resort to the theory developed
by Auroux and Katzarkov  (\cite{a-k}),
secondly, we prove  more general results than the ones
contained in \cite{cat02}.
This is needed
in order to show  that
the Manetti examples are canonically symplectomorphic.

The examples of \cite{man4}
are in fact pairs of surfaces of
general type $S, S'$ which are not deformation equivalent, but which admit
respective degenerations to  normal surfaces $X, X'$ with singularities,
which in turn are deformation equivalent to each other via an equisingular
deformation.

Moreover, each  degeneration yields
   a $\Q$-Gorenstein smoothing of the singularities.
Manetti used  a result of Bonahon on the group of diffeomorphisms
of lens spaces in order to show that the surfaces $S, S'$
are diffeomorphic to each other.
Here,  the canonical symplectomorphism of Manetti surfaces
follows from the  following

\begin{teo}\label{families}
Let $  \mathcal {X}\subset \PP ^N \times \Delta$ and
$  \mathcal {X}' \subset \PP ^N \times \Delta'$ be two flat families
of normal surfaces over the disc of radius 2 in $\C$.

Denote by $ \pi :  \mathcal {X} \to  \Delta$ and
by  $ \pi ' :  \mathcal {X}' \to  \Delta$ the respective projections
and make the following assumptions on the respective fibres of $\pi , \pi'$:

1) the central fibres
$X_0 $ and $ X_0'$ are surfaces with cyclic quotient singularities
and the two flat families yield $\Q$-Gorenstein smoothings of them,

2) the other fibres $X_t$, $X_t'$, for $ t, t' \neq 0$ are smooth.

Assume moreover that

3) the central fibres
$X_0 $ and $ X_0'$ are projectively equivalent to
respective fibres ($X_0 \cong Y_0$ and $ X_0' \cong Y_1$)
   of an equisingular
projective family $  \mathcal {Y}\subset \PP^N \times \Delta$
of surfaces.

Set $X := X_1$, $X' := X_1'$: then

a) $X$ and $X'$ are diffeomorphic

b) if  $FS $ denotes the symplectic form inherited from the Fubini-Study
K\"ahler metric on $\PP^N$, then the symplectic manifolds
$(X,FS)$ and $(X',FS)$ are symplectomorphic.

The same conclusion holds if  hypothesis 1) is replaced by

1') The singularities of $X_0 $,
resp. $ X'_0$, admit a unique smoothing component.
\end{teo}

An important step in the proof is furnished by the following

\begin{teo}\label{glueing}
Let $X_0 \subset \PP^N$ be a projective variety with isolated
    singularities admitting a smoothing component.

Assume that, for each singular point $x_h \in X$,
we choose a smoothing component $T_{j(h)}$ in the basis of the
semiuniversal deformation of the germ $(X, x_h)$.
Then (obtaining different results for each such choice) $X$ can be
approximated (in the Hausdorff metric) by symplectic submanifolds
$W_t$ of $ \PP^N$, which are diffeomorphic to the glueing
of the 'exterior' of $X_0$ (the complement to the union $B = \cup_h B_{h}$ of
     suitable (Milnor) balls around the singular points) with
the  Milnor fibres $\sM_h$ , glued along the singularity links $\sK_{h,0}$.
\end{teo}

A pictorial view of the proof is contained in the following Figure.


A corollary of the  above theorems is the following

\begin{teo}\label{canonical}A minimal surface of general type
$S$ has a canonical symplectic
   structure, unique up to symplectomorphism, and stable by deformation,
such that the class of
the symplectic form is the class of the canonical sheaf
$ \omega_S =  \Omega^2_S = \hol_S (K_S).$
The same result holds for any projective smooth variety with
ample canonical bundle.

\end{teo}

\begin{oss}

Theorem \ref{families} holds
   more generally for varieties of higher dimension with single smoothing
isolated singularities,
or under the assumptions

i)  $X_0 = X'_0$

ii) for each singular point $ x_0$ of $X_0$,
the two smoothings $  \mathcal {X},   \mathcal {X}'$,
correspond to paths in the same irreducible
component of $Def (X_0, x_0) $.
\end{oss}

As  already mentioned,  Manetti surfaces are Abelian
coverings  of rational surfaces with group $( \Z /2)^m$, leading to
the situation of  Theorem 1.2. Hence

\begin{teo}\label{manetti}
Manetti's surfaces (Section 6 in \cite{man4})
provide examples of surfaces of general type
which are not deformation equivalent, but, endowed with their
canonical symplectic structures, are symplectomorphic.
\end{teo}

\section{Preliminaries}

Let us first recall the well known Theorems of Ehresmann (\cite{ehre})
and Moser (\cite{moser})

\begin{teo}\label{e-m}(Ehresmann + Moser)
Let $ \pi :  \mathcal {X}\rightarrow T  $ be a proper submersion of
differentiable manifolds with $T$ connected, and assume that we have
a differentiable $2$-form $\omega$ on $\mathcal {X}$ with the property
that

(*) $\forall t \in T$ $ \omega_t : = \omega |_{X_t}$ yields
a symplectic structure on $X_t$ whose class in $H^2 (X_t, \R)$
   is locally constant on $T$ (e.g., if it lies on $H^2 (X_t, \Z)$).

Then the symplectic manifolds $(X_t, \omega_t)$ are all symplectomorphic.
\end{teo}

There is  an easy  variant with boundary of
Ehresmann's theorem

\begin{lem}\label{e-v}
Let $ \pi :  \mathcal {M} \rightarrow T  $ be a proper submersion of
differentiable manifolds with boundary, such that $T$ is
a ball in $\R^n$, and
assume that we are given a fixed trivialization  $\psi$ of
   a  closed family $  \mathcal {N} \rightarrow T  $ of submanifolds with
boundary. Then we can find a
trivialization of $ \pi :  \mathcal {M} \rightarrow T  $ which
induces the given trivialization $\psi$.
\end{lem}

\begin{proof}
It suffices to take on $\mathcal {M}$ a Riemannian metric
where the sections $ \psi (p,T) $, for $ p \in  \mathcal {N}$,
are orthogonal to the fibres of $\pi$. Then we use the customary proof
of Ehresmann's theorem, integrating  liftings orthogonal to the fibres of
   standard vector fields on $T$.
\end{proof}

The variant \ref{e-v} of Ehresmann's theorem will now be first
applied to the
Milnor fibres of smoothings of isolated singularities.

Let $ (X,x_0)$ be the germ of an isolated singularity of a complex space,
which is pure dimensional
of dimension $ n = dim_{\C}X$,
assume $x_0 = 0 \in X \subset \C^{n+m}$, and consider
the ball $ \overline{B ( x_0, \delta)}$ with centre
the origin and radius $\delta $.
Then, for all $ 0 < \delta << 1$, the intersection $\sK_0 : =  X \cap S ( x_0,
\delta)$, called the {\bf link} of the singularity, is a smooth manifold of
real dimension $ 2n-1$.

Consider the semiuniversal deformation $ \pi : (\X, X_0, x_0) \ra
(\C^{n+m},0) \times
(T, t_0)$ of $X$ and
the family of singularity links
$\sK : =  \X \cap  (S ( x_0, \delta) \times
(T, t_0))$. By a uniform continuity argument it follows that
$\sK \to T$ is a trivial bundle if we restrict $T$ suitably around
the origin $ t_0$
(it is a differentiably trivial fibre bundle in the sense of
stratified spaces, cf. \cite{math}).

We can now recall the concept of Milnor fibres of $ (X,x_0)$.

\begin{df}\label{Milnor}
Let $(T, t_0)$ be the basis of the semiuniversal deformation of a germ
of isolated singularity  $ (X,x_0)$, and let $T = T_1 \cup \dots \cup T_r$
be the decomposition of $T$ into irreducible components.
$T_j$ is said to be a smoothing component if there is a $ t \in T_j$
such that the
corresponding fibre $X_t$ is smooth.  If $T_j$ is a smoothing component,
then the corresponding Milnor fibre is the intersection
of the ball $ \overline{B ( x_0, \delta)}$ with the fibre $X_t$, for
$ t \in T_j$, $ |t| < \eta << \de << 1$.
\end{df}

Whereas the singularity links form a trivial bundle, the Milnor fibres
form only a differentiable bundle of manifolds with boundary
over the open set
$ T^0_j : = \{ t \in T_j ,  |t| < \eta | X_t {\rm \ is \ smooth} \}.$

Since however $T_j$ is irreducible, $ T^0_j$ is connected, and the Milnor fibre
is unique up to smooth isotopy, in particular up to diffeomorphism.

   Let us now recall some known facts on the
   class of singularities given by the (cyclic)
quotient singularities admitting a $\Q$-Gorenstein smoothing
(cf. \cite{man4}, Section 1, pages 34-35, or the original sources
\cite{k-sb},\cite{man0},\cite{l-w}).

The simplest way to describe the singularities
    $${\bf Cyclic \ quotient \ singularity} \ \frac{1}{d n^2} (1,dna-1)
   = A_{dn-1} / \mu_n$$
is to view them on the one side as quotients of $\C^2$ by a cyclic
group of order $d n^2$ acting with the indicated characters $(1,dna-1)$,
or on the other side as quotients  of the rational double point $ A_{dn-1}  $
of equation $ uv - z^{dn} = 0$ by the action of the group $\mu_n$ of
n-roots of unity acting in the following way:

$$\xi \in \mu_n  \ {\bf acts \  by \ :} \
   (u,v,z) \rightarrow  ( \xi u , \xi^{-1} v,\xi^{a} z).$$

This quotient action gives rise to a quotient family
$ \mathcal{X} \rightarrow \C^d$, where

$ \mathcal{X}=
   \mathcal{Y}/ \mu_n$ , $ \mathcal{Y}$ is the hypersurface in
$\C^3 \times \C^d$ of equation
$ uv - z^{dn} = \Sigma_{k=0}^{d-1} t_k z ^{kn}$ and the action of
$\mu_n$ is extended trivially on the factor $\C^d$.

The heart of the construction is that $\mathcal{Y}$, being
   a hypersurface,
is Gorenstein (this means that the canonical sheaf
   $\omega_{\mathcal{Y}}$ is
invertible), whence such a quotient $ \mathcal{X}=
   \mathcal{Y}/ \mu_n$, by an action which is unramified in codimension $1$,
is (by definition) $\Q$-Gorenstein.

These smoothings were considered by Koll\'ar and Shepherd Barron
   (\cite{k-sb}, 3.7-3.8-3.9,  cf. also   \cite{man0}),
who pointed out their relevance in the theory
of compactifications of moduli spaces of surfaces,
and showed that, conversely, any $\Q$-Gorenstein smoothing
of a quotient singularity is induced by the above family (which has a
smooth base, $\C^d$).

Riemenschneider (\cite{riem})  had earlier shown that,
   for the cyclic quotient
singularity $\frac14 (1,1)$, the basis of the
semiuniversal deformation   consists of two smooth components
   intersecting transversally, each one yielding a smoothing, but
only one admitting a simultaneous resolution, and only the other
yielding   smoothings with $\Q$-Gorenstein total space.

\section{Proof of the Theorems }
\label{second}
\begin{proof}{\bf (of Theorem \ref{families}) }

Applying Theorem \ref{e-m} to $T: = \Delta - \{0\}$,
and to the restrictions of the two given families $\mathcal {X}$,
$\mathcal {X}'$, we can for both statements replace $X$ by any
$X_t$ with $t \neq 0$ sufficiently small, and similarly replace $X'$ by
any $X'_{t'}$ with $t' \neq 0$.

In other words, assuming $X_0 , X'_0 \subset \PP^N$, we may assume that
$X$ and $X'$ are  very near to $X_0$, respectively $X'_0$.

Since the family $\mathcal {Y}$ is equisingular, to each singular
point $x_0 \in Sing(X_0)$ corresponds a unique singular point
$x'_0 \in Sing(X'_0)$(indeed $X_0, X'_0$ are homeomorphic by
a homeomorphism carrying $x_0$ to $x'_0$).

For each $x_0 \in Sing(X_0)$, $\pi$  induces a germ of
holomorphic mapping
$F_{x_0} : \Delta \rightarrow \mathcal{D}_{x_0} \subset Def (X_0, x_0)$,
where $\mathcal{D}_{x_0}$ is the chosen smoothing component;
respectively,  for $x'_0 \in Sing(X_0)$ the corresponding singular point,
   $\pi'$ induces
a germ  $F'_{x'_0}$.

Let $\mathcal{Z}_{x_0} \subset \mathcal{D}_{x_0} \times \PP^N$ be
given by the restriction of the semiuniversal deformation of the germ
$(X_0, x_0)$,
and, for each $ t \in \Delta$,  consider the corresponding
singular point $y_0(t) \in Y_t$ (thus $y_0(0)= x_0, y_0(1) = x'_0 $),
and the semiuniversal deformation of the corresponding germ
$(Y_t , y_0(t) )$.

We obtain in this way a family of pairs of germs,
$$\Y'_0 \subset \mathcal{Z}_{0} \subset \mathcal{D}_{0}
\times \PP^N \ra \mathcal{D}_{0}
   \subset \Delta \times \C^m,$$
where $\Y'_0$ is the family of germs $(Y_t , y_0(t) )$, induced by $\Y$.

Our assumption, that each $\mathcal{D}_{y_0(t)}$ is irreducible,
implies immediately the irreducibility
of $\mathcal{D}_{0}$.

For each $ 0 < \epsilon << 1 ,0 < \eta << 1$
we consider the family of Milnor
links

$$ \K_{\epsilon, \eta} : =  \cup_t \mathcal{K}_{\epsilon, \eta} (t):=
\cup_t [ \mathcal{Z}_{0} \cap
   (\{ t\} \times B(0, \epsilon ) \times   S (y_0 (t), \eta ) )]$$

where $B(0, \epsilon )$ is the ball of radius $\epsilon$ and centre
the point $0 \in \mathcal{D}_{y_0(t)}$ corresponding to $(Y_t , y_0(t) )$,
while $S (y_0 (t), \eta ) $ is the sphere in $\PP^N$ with centre $y_0 (t)$ and
radius $\eta$ in the Fubini Study metric.

We already remarked that, for $\eta << 1$ and $ \epsilon << \eta $,
the family $ \mathcal{K}_{\epsilon, \eta} \rightarrow
   ((\De \times B(0, \epsilon )) \cap \mathcal{D}_{0}) $ is differentially
trivial (either in the sense of stratified sets, cf. \cite{math},
or, as suffices to us, in the weaker sense that
when we pull it back through a differentiable
map $\Delta \rightarrow  (B(0, \epsilon ) \cap \mathcal{D}_{x_0}) $
we get a differentiable product).

{\bf Proof of a):} we apply lemma \ref{e-v} several times:

\begin{itemize}
\item
i) first  we apply it in order
to thicken the trivialization of the Milnor links to a closed tubular
neighbourhood
in the family $\ZZZ_0$,
\item
ii-a) then we apply it in order to get a compatible trivialization of the
family of exteriors in $Y_t$ of the balls $B ( y_0(t), \eta /2)$
\item
ii-b) then we apply it to the restriction of the families
$\mathcal {X} \rightarrow \Delta  $, $\mathcal {X}' \rightarrow \Delta  $,
to a ball of radius  $\delta$ where  $\delta$  is so chosen that
$F_{x_0}  (\{ t | \ |t| < \delta \} ) \subset B(0, \epsilon / 2)$
(resp. for $F'_{x_0}$), and to the exterior of the
balls $B ( x_0, \eta /2)$,resp. $B ( x'_0, \eta /2)$,
   so that we  get compatible trivializations
   of the exterior in $\X$
to the balls $B ( x_0, \eta /2)$, resp. in $\X'$
to the balls $B ( x'_0, \eta /2)$
\item
iii) we finally use our assumptions that the images of $F'_{x_0}$, resp.
$F_{x_0}$ land in  $\mathcal{D}_{0}$ which is irreducible:
it follows that there is a holomorphic mapping $ G: \Delta
\rightarrow \mathcal{D}_{0}$ whose image contains the
two points $F_{x_0} (t_0)$, $F'_{x'_0}(t'_0)$ and is contained in
$ ( \De \times B(0, \epsilon / 2) ) \cap \Sigma $ ( $\Sigma$ being as before
   the smoothing locus).

We consider then the pull back to $\Delta$ under $G$ of the family of
   closed Milnor fibres

$$ \mathcal{M}_{\epsilon, \eta} :=
   \cup_t \mathcal{M}_{\epsilon, \eta} (t):=
\cup_t [ \mathcal{Z}_{0} \cap
   (\{ t\} \times B(0, \epsilon ) \times  \overline{ B (y_0 (t), \eta )} )].$$
To this family we apply again 2.2, in order to obtain a trivialization
of the family of Milnor fibres which extends the given trivialization on
the family of (closed) tubular neighbourhoods of the Milnor links.
\end{itemize}

We are now done, since we obtain the desired diffeomorphism between
$X$ and $X'$ by glueing together (in the intersection with
$B (x_0, \eta ) - \overline{B (x_0, \eta /2 )}$,
resp. with $B (x'_0, \eta ) - \overline{B (x'_0, \eta /2 )}$ )
   the two diffeomorphisms
provided by the restrictions of the respective trivializations
ii) (to the intersection of the complement to
$\overline{B (x_0, \eta /2 )}$, resp. $\overline{B (x'_0, \eta /2 )}$)
and iii) ( to the intersection
with  $B (x_0, \eta )$, resp.$B (x'_0, \eta )$ ):
they glue because they both extend
   the trivialization i).

{\bf Step I in the proof of b. } We want to use the previous construction and
considerations in order
to construct a family of differentiable embeddings of the same
differentiable manifold $V \cong X \cong X'$ into $\PP^N$,
which includes the embeddings $X$ and $X'$. We shall later show in step II
that we can manage  that
every fibre inherits a symplectic structure from the Fubini-Study
form and we can then finally apply Moser's theorem.

First of all, observe that $V$ is obtained from $X$ by writing $X$ as the
glueing of two manifolds with boundary, namely the union $M$
of the Milnor fibres,
   and the 'exterior' of $X$, i.e., the closure $\E$  of the complement
$X \setminus B$ ($B=$ union of the Milnor balls
$B(x_0, \eta)$), which both have as boundary the union $\K$
of the Milnor links.

We consider now the product  $ V \times [-1, 1]$ and
we first map it to $\PP^N \times [-1, 1]$ via the product of
a  piecewise differentiable
map  $\psi$ and the second projection. Later on we shall approximate
$\psi$ by a differentiable map $\phi$ which is an embedding
when restricted to each fibre.

The key idea to construct $\psi$ is to make a little bit longer the neck
around the Milnor link, and to use the trivialization of
the family of Milnor links and of Milnor fibres.

We define $\psi_s$, for $ -1 \leq s \leq - \frac12$ on $\E$ by using
a path $\tau(s)$ in $\D_{x_0}$ from the point $t$ corresponding to $X$,
and reaching the origin ($X_0$) for $ s = - \frac12$: thus  $t := \tau(-1)$,
and $ 0 := \tau (-\frac12)$. We map then the exterior $\E$ to the
exterior of $X_{ F_{x_0}(\tau (s))} $, we do a completely similar operation
for $ \frac12 \leq s \leq 1$ (thus, for $s=1$ $\E$ maps to the exterior
of $X'$).

Instead, for $ - \frac12 \leq s \leq  \frac12 $, we use a path
$\nu(s)$ in the parameter space for the family $\Y$, and
map the exterior $\E$ to the exterior of $Y_{\nu(s)}$.

For the interiors, i.e., the Milnor fibres minus a collar $\CCC$, we
send them to the Milnor fibres corresponding to a path
$\sigma (s)$ in $\D_0$ connecting
the points $F_{x_0}(t), F'_{x'_0}(t')$ corresponding to $X$,
respectively to $X'$. We require of course that the first
projection of $\sigma(s)$ onto $\De$ equals $\nu(s)$ for
$ - \frac12 \leq s \leq  \frac12 $, while for
$ -1 \leq s \leq - \frac12$ the first coordinate of $\sigma (s)$
equals $0$, similarly for $ \frac12 \leq s \leq 1$ it equals $1$.

Finally, we use the collar and the trivialization of the Milnor links to
join the boundary of the Milnor fibres with the boundary of the exterior.

We have now, for each $s$, a map $\psi_s$ which is a differentiable
embedding in the
exterior $ \forall s$, while it is an embedding for $ s= -1, +1$.

Without loss of generality we may assume that $ N \geq 4$,
else there is nothing to prove,
and we indeed may assume $ N \geq 5$ by choosing if $ N=4$
the standard embedding $ \PP^4 \ra \PP^5$.

We obtain the desired family of embeddings $\phi_t$ by applying
the following variation of Whitney's embedding theorem

\begin{lem}
Let W be a differentiable manifold with boundary of dimension $k$, and let
$\psi : W \ra \R^n$ be a continuos map which is
a differentiable embedding around the boundary $\partial W$.

If $ n \geq 2k+1$, then $\psi$ can be approximated by an embedding
$\phi$ which coincides with $\psi$ in a neighbourhood
of $\partial W$.
\end{lem}
{\em Proof of the Lemma}
First of all, by a standard convolution process,
we can approximate $\psi$ by  a differentiable map
$ \varphi : W \ra \R^n$ which coincides with $\psi$ in a neighbourhood
of $\partial W$.

After that, we construct a  differentiable map $f : W \ra \R^m$
which equals zero in a neighbourhood
of $\partial W$, and is such that
$ \varphi \oplus f: W \ra \R^n \oplus \R^m$ is an embedding.

After that, we construct $\phi$ composing $ \varphi \oplus f$
via a sufficiently general
projection $\R^n \oplus \R^m \ra \R^n $ close to the first
projection of the direct sum.

\end{proof}

We are not yet finished with the proof of b), because we have that the
pull back $\omega$ of the Fubini Study form is non degenerate
on the exterior of the Milnor balls, and in the interior of the
Milnor fibres, but it could be degenerate in a neighbourhood
of the collar of the Milnor link $\K$.

The procedure to obtain nondegeneracy everywhere will be now given
while proving  Theorem
\ref{glueing}.

\begin{proof}{\bf (of Theorem \ref{glueing}) }
Without loss of generality, and to simplify our notation, assume that
$X_0$ has only one singular point $x_0$, and let $B: = B(x_0, \eta)$
be a Milnor ball around the singularity, let moreover $\D := \D_{x_0} \subset
Def (X_0, x_0)$, and, for $t \in \D \cap B(0, \epsilon)$ we consider the
Milnor fibre $\MM_{\epsilon,\eta}(t)$, whereas we have the two Milnor links
   $$\K_0 : = X_0 \cap S(x_0, \eta)  \ {\rm and } \ \K_t : =
   \ZZZ_t  \cap S(x_0, \eta - \de) $$ .

We can consider the Milnor collars $\CCC_0(\de) := X_0 \cap
(\overline{B(x_0, \eta)} \setminus B(x_0, \eta - \de))$,
and $\CCC_t(\de) := \ZZZ_t \cap
(\overline{B(x_0, \eta)} \setminus B(x_0, \eta - \de))$.

We restrict now $t$ to vary in a holomorphic disc $\De$ mapping
to a path through the origin and intersecting the complement of
the smoothing locus $\Sigma$ only at the origin.

Now, with this restriction, the Milnor collars fill up a complex
submanifold of dimension $ dim X_0 +1 := n + 1$, and
we 'join the two Milnor links'  by a differentiable
embedding  of the  abstract Milnor collar (i.e., $\CCC_0(\de)$
as a differentiable manifold maps in such a way that its boundary
maps onto $\K_0  \cup \K_t$).

For $\epsilon < < \de $ the tangent spaces to the image of the
abstract Milnor collar can be made very close to the tangent spaces
of the Milnor collars $\MM_{\epsilon,\eta}(t)$, and we can conclude
the proof via the following well known
lemma which lies at the heart of Donaldson's
work (\cite{don6})

\begin{lem}
Let $W \subset \PP^N$ be a differentiable submanifold of
even dimension $2n$, and assume that the tangent space of $W$
is close to be complex in the  sense that for each tangent vector
$v$ to $W$ there is another tangent vector $v'$ such that
$$  J v = v' + u, | u | <  | v| .$$

Then the restriction to $W$ of the Fubini Study Form $ \omega_{FS}$
makes $W$ a symplectic submanifold of $\PP^N$.
\end{lem}

\Proof
Let $A$ be the symplectic form on projective space, so that for each
tangent vector $v$ we have:

$  | v|^2 = A ( v, J v) = A (v, v') + A (v, u) .$

Since $ | A (v, u) | < |v|^2 $, $ A (v, v') \neq 0$ and
$A$ restricts to a nondegenerate form.
\end{proof}

\QED for Theorem \ref{glueing}

\begin{oss}
By Moser's theorem the symplectic manifolds $W_t$ are symplectomorphic to
each other.
\end{oss}

{\bf Step II in the proof of b. }

We apply the same method used as in the proof of Theorem \ref{glueing},
i.e., when we apply the 'relative'  Whitney theorem to glue
the Milnor fibres corresponding to points $\sigma(s)$,
   we choose $|\sigma(s)| < \epsilon$,  $| F_{x_0}(t)| < \epsilon$,
$| F'_{x'_0}(t')| < \epsilon$
and , using the compactness
of the interval $[- \frac12 , + \frac12 ]$, there exists
$ \epsilon$ so small such that the tangent spaces of the submanifolds
$W_s$ are close to be complex, hence the $W_s$ are symplectic
submanifolds.

Then Moser's theorem gives the symplectomorphism of $X$ with $X'$.

\QED for Theorem \ref{families}

\begin{proof}{\bf (of Theorem \ref{canonical}) }
Let $S$ be the minimal model of a surface of general type.

We prove the assertion first in any dimension, but in
the special case where we have a smooth projective variety $V$ whose  canonical
divisor $ K_V$ is ample.

In fact, let $m$ be such that $m K_V$ is very ample
(any $m\geq 5$ does by Bombieri's theorem, cf. \cite{bom} in the case
of surfaces, for higher dimension we can use Matsusaka's big theorem,
   cf. \cite{siu} for an effective version) thus
the $m$-th pluricanonical map $ \phi_m := \phi_{|mK_V|}$
is an embedding of $V$ in a projective space $\PP^{P_m-1}$,
where $P_m := dim H^0(\hol_S (m K_V) )$.

We define then  $\omega_m$ as follows:  $\omega_m := \frac{1}{m} \phi_m
   ^* (FS)$ (where $FS$
is the Fubini-Study form $\frac{1}{2 \pi i} \partial
\overline{\partial} log |z|^2$), whence $\omega_m$ yields
a symplectic form as desired.

One needs to show that the symplectomorphism class of $(V, \omega_m)$
is independent of $m$. To this purpose, suppose that $n$ is also such
that $ \phi_n$  yields an embedding of $V$: the same holds also for
$nm$, whence it suffices to show that $(V, \omega_m)$ and
$(V, \omega_{mn})$ are symplectomorphic.

To this purpose we use first the well known and easy fact
that the pull back of the Fubini-Study form under the $n$-th
Veronese embedding $v_n$ equals the $n$-th multiple of the Fubini-Study
form. Second,  since $v_n \circ \phi_m$ is a linear projection
of $\phi_{nm}$, by Moser's Theorem follows the desired symplectomorphism.

Assume now that $S$ is a minimal surface of general type and
that $ K_S$ is not ample: then for any $m\geq 5$ (by Bombieri's
cited theorem) $ \phi_m$  yields an embedding of the canonical model
$X$ of $S$, which is obtained by contracting the finite number of
smooth rational curves with selfintersection number $= -2$ to a finite
number of Rational Double Point singularities. For these, the
base of the semiuniversal deformation is smooth and yields a
   smoothing of the singularity.

By Tjurina's theorem (cf. \cite{tju}), $S$
is diffeomorphic to any smoothing $S'$ of $X$:
however we have to be careful because there are many examples
(cf. e.g. \cite{cat5}) where
$X$ does not admit any global smoothing.

Since however there are local smoothings,  Tjurina's theorem
tells us that $S$ is obtained glueing the exterior
$X \setminus B$ ($B$ being the union of Milnor balls of radius $\eta$ around
the singular points of $X$) together  with the respective Milnor fibres.

Argueing as in theorem \ref{glueing} we find a symplectic submanifold $W$
of projective space which is diffeomorphic to $S$, and by the previous remark
$W$ is unique up to symplectomorphism. Clearly moreover, if $X$ admits
a smoothing, we can then take $S'$ sufficiently close to $X$
   as our approximation $W$.
Then $S'$ is a surface with ample canonical bundle, and, as we have
seen, the symplectic structure induced by (a submultiple of) the Fubini Study
form is the canonical symplectic structure.

The stability by deformation is again a consequence of Moser's theorem.

\end{proof}

\begin{proof}{\bf  (of Theorem \ref{manetti})}
In \cite{man4} Manetti constructs examples of surfaces $S$, $S'$
of general type which are not deformation equivalent, yet with
the property that
there are flat families of normal surfaces
$  \mathcal {X}\subset \PP^N \times \Delta$ and
$  \mathcal {X}' \subset \PP^N \times \Delta'$

1) yielding a $\Q$-Gorenstein smoothings of
   the central fibres $X_0 , X_0'$,

and such that

2) the  fibres $X_t$, $X_t'$, for $ t, t' \neq 0$ are smooth,
and the canonical divisor of each fibre is ample

3) there are $t_0$, $t'_0$ with
$S \cong X_{t_0}$, $S' \cong X'_{t'_0}$

4) there exists an equisingular family $\Y \ra \De$ whose fibres have
indeed only singularities of type $\frac14 (1,1)$, and such that
$ Y_0 \cong X_0$, $ Y_1 \cong X'_0$.

There exists therefore a positive integer $m$ such that for
each $X_t$ and $X'_{t'}$ the $m$-th multiple of the canonical
(Weil-)divisor is Cartier and very ample, and therefore
the relative $m$-pluricanonical maps yield three new projective families
to which 1.2 applies.

By 1.2 and 1.3 it follows that $S$ and $S'$, endowed with their canonical
symplectic structure, are symplectomorphic.

\end{proof}

\bigskip

{\bf Acknowledgements.}

I would like to thank the referee of \cite{cat4}
for prompting me to write
down this note, and Marco Manetti for observing
that the result proven in \cite{cat02} was not
sufficient for the application to Theorem \ref{manetti},
since his surfaces do not admit a direct degeneration to the same
normal surface.

The research of the  author was performed in the realm  of the
   SCHWERPUNKT "Globale Methoden in der komplexen Geometrie",
and of the EAGER EEC Project.

Thanks to the IHES for hospitality in the month of march 2006, where the
final version of the article was written; the visit was
supported by contract nr. RITA-CT-2004-505493,
project acronym: IHES Euro-Programme 2.

\bigskip

\begin{footnotesize}
\noindent
{\bf Notes.}
a) Marco Manetti  informed me that he was aware of a result
similar to part a) of Theorem 1.2.

b) similar ideas to the ones of the present paper, but with different proofs
and especially with different applications appear in \cite{sty02} and
\cite{st03}.
\end{footnotesize}

\vfill

\noindent
{\bf Author's address:}

\bigskip

\noindent
Prof. Fabrizio Catanese\\
Lehrstuhl Mathematik VIII\\
Universit\"at Bayreuth, NWII\\
   D-95440 Bayreuth, Germany

e-mail: Fabrizio.Catanese@uni-bayreuth.de

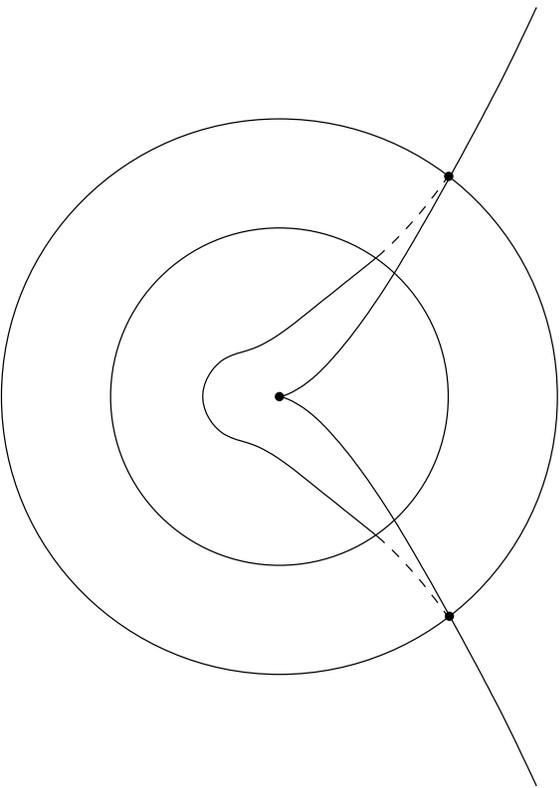
\end{document}